\documentclass[conference]{IEEEtran}
\usepackage{amsmath, amssymb}

\newtheorem{theorem}{Theorem}[section]

\newtheorem{lemma}[theorem]{Lemma}

\newtheorem{remark}[theorem]{Remark}

\def \< {\langle}
\def \> {\rangle}
\def \id {{\it id}}

\def \R {\mathbb{R}}

\def \C {\mathbb{C}}
\def \E {\mathbb{E}}
\def \supp {{\rm supp}}
\def \Ker {{\rm Ker}}
\def \sign {{\rm sign}}
\def \polylog {{\rm polylog}}
\def \conv {{\rm conv}}
\def \etc {,\ldots,}
\def \e {\varepsilon}
\def \d {\delta}
\def \w {\omega}

\newcommand{\pr}[2]{\langle {#1} , {#2} \rangle}
\newcommand{\norm}[1]{\Big\| #1 \Big\|}

\begin{document}

\title{Sparse reconstruction by convex relaxation:\\
Fourier and Gaussian measurements}

\author{\authorblockN{Mark Rudelson}
\authorblockA{Department of Mathematics \\
   University of Missouri, Columbia \\
   Columbia, Missouri 65211 \\
   Email: rudelson@math.missouri.edu}
\and
\authorblockN{Roman Vershynin}
\authorblockA{Department of Mathematics \\
  University of California, Davis \\
  Davis, California 95616 \\
  Email: vershynin@math.ucdavis.edu}
}
\maketitle

\begin{abstract}
This paper proves best known guarantees for exact reconstruction
of a sparse signal $f$ from few non-adaptive universal linear
measurements.
We consider Fourier measurements (random sample of frequencies
of $f$) and random Gaussian measurements. The method for reconstruction
that has recently gained momentum in the Sparse Approximation Theory
is to relax this highly non-convex problem to a
convex problem, and then solve it as a linear program.
What are best guarantees for the reconstruction problem
to be equivalent to its convex relaxation is an open question.
Recent work shows that the number of measurements $k(r,n)$
needed to exactly reconstruct any $r$-sparse signal $f$ of length $n$
from its linear measurements with convex relaxation
is usually $O(r \, \polylog(n))$.
However, known guarantees involve huge constants, in spite of
very good performance of the algorithms in practice.
In attempt to reconcile theory with practice, we prove the first
guarantees for universal measurements (i.e. which work for all 
sparse functions) with reasonable constants. 
For Gaussian measurements,
$k(r,n) \lesssim  11.7 \, r \big[ 1.5 + \log(n/r) \big]$,
which is optimal up to constants.
For Fourier measurements, we prove the best known bound
$k(r,n) =  O(r \log(n) \cdot \log^2(r) \log(r \log n) )$,
which is optimal within the $\log \log n$ and $\log^3 r$ factors.
Our arguments are based on the technique of Geometric Functional 
Analysis and Probability in Banach spaces.
\end{abstract}

\IEEEpeerreviewmaketitle

\section{Introduction}
During the last two years, the Sparse Approximation Theory benefited
from a rapid development of methods based on the Linear Programming. The idea was to
relax a sparse recovery problem to a convex optimization problem. The convex problem
can be further be rendered as a linear program, and analyzed with all available
methods of Linear Programming.

Convex relaxation of sparse recovery problems can be traced back
in its rudimentary form to mid-seventies; references to its early
history can be found in \cite{T 04c}. With the development of fast
methods of Linear Programming in the eighties, the idea of convex
relaxation became truly promising. It was put forward most
enthusiastically and successfully by Donoho and his collaborators
since the late eighties, starting from the seminal paper \cite{DS}
(see Theorem~8, attributed there to Logan, and Theorem~9).
There is extensive work being carried out, both in theory and in
practice, based on the convex relaxation
\cite{CDS, DH, EB, FN, DE, GN, T 04a, T 04b, T 04c, DET, D 04a, D 04b, 
DT 04a, CRT, CR, CT, CT 05, RV, CRTV, CT 06, MPT}.

To have theoretical guarantees for the convex relaxation method,
one needs to show that {\em the sparse approximation problem is
equivalent to its convex relaxation}. Proving this presents a
mathematical challenge. Known theoretical guarantees work only for
random measurements (e.g. random Gaussian and Fourier
measurements). Even when there is a theoretical guarantee, it
involves intractable or very large constants, far worse than in
the observed practical performances.

In this paper, we substantially improve best known theoretical guarantees for random
Gaussian and Fourier (and non-harmonic Fourier) measurements.
For the first time, we are able to prove guarantees
with reasonable constants (although only for Gaussian measurements).
Our proofs are based on methods of Geometric Functional Analysis,
Such methods were recently successfully used for related problems
\cite{RV}, \cite{MPT}.
As a result, our proofs are reasonably short (and hopefully, transparent).

In Section~\ref{s:relax}, we state the sparse reconstruction problem and
describe the convex relaxation method. A guarantee of its
correctness is a very general {\em restricted isometry condition}
on the measurement ensemble, due to Candes and Tao (\cite{CT 05}, see \cite{CRTV}).
Under this condition, the reconstruction problem with respect to
these measurements is equivalent to its convex relaxation. In
Sections~\ref{s:Fourier} and \ref{s:Gauss},
we improve best known guarantees for the sparse
reconstruction from random Fourier (and non-harmonic Fourier)
measurements and Gaussian measurements
(Theorem~\ref{Fourier rec} and \ref{Gaussian rec} respectively).

\section{The Sparse Reconstruction Problem and its Convex Relaxation} \label{s:relax}

We want to reconstruct an unknown signal $f \in \C^n$ from linear measurements
$\Phi f \in \C^k$, where $\Phi$ is some known $k \times n$ matrix,
called the {\em measurement matrix}. In the interesting case $k < n$, the
problem is underdetermined, and we are interested in the sparsest solution.
We can state this as the optimization problem
\begin{equation}                    \label{non-convex}
  \text{minimize }  \|f^*\|_0 \text{ subject to } \Phi f^* = \Phi f,
\end{equation}
where $\|f\|_0 = |\supp f|$ is the number of nonzero coefficients of $f$.
This problem is highly non-convex. So we will consider its {\em convex relaxation}:
\begin{equation}                    \label{convex}
  \text{minimize }  \|f^*\|_1 \text{ subject to } \Phi f^* = \Phi f,
\end{equation}
where $\|f\|_p$ denotes the $\ell_p$ norm throughout this paper,
$(\sum_{i=1}^n |f_i|^p)^{1/p}$.
Problem \eqref{convex} can be classically reformulated as the {\em linear program}
$$
\text{minimize }  \sum_{i=1}^n t_i
\text{ subject to } -t \le f^* \le t, \ \Phi f^* = \Phi f,
$$
which can be efficiently solved using general or special methods
of Linear Programming. Then the main question is:
\begin{quote}
  {\em Under what conditions on $\Phi$ are problems \eqref{non-convex}
  and \eqref{convex} equivalent?}
\end{quote}
In this paper, we will be interested in the {\em exact reconstruction}, i.e.
we expect that the solutions to \eqref{non-convex} and \eqref{convex} are
equal to each other and to $f$. Results for approximate reconstruction
can be derived as consequences, see \cite{CT}.

For exact reconstruction to be possible at all, one has to assume that the
signal $f$ is $r$-sparse, that is $\supp(f) \le r$, and that the number
of measurements $k = k(r,n)$ has to be at least twice the sparsity $r$.
Our goal will be to find sufficient conditions (guarantees) for the exact
reconstruction. The number of measurements $k(r,n)$ should be kept as small as
possible. Intuitively, the number of measurements should be of the order of
$r$, which is the `true' dimension of $f$, rather than the nominal
dimension $n$.

Various results that appeared over the last two years demonstrate
that many natural measurement matrices $\Phi$ yield exact
reconstruction, with the number of measurements $k(r,n) = O(r
\cdot \polylog(n))$, see \cite{CRT, CT, CT 05, RV}.
In Sections~\ref{s:Fourier} and \ref{s:Gauss},
we improve best known estimates on $k$ for Fourier (and, more
generally, nonharmonic Fourier) and Gaussian matrices respectively.

A general sufficient condition for exact reconstruction is the
{\em restricted isometry condition} on $\Phi$, due to Candes and
Tao (\cite{CT 05}, see \cite{CRTV}).
It roughly says that the matrix $\Phi$ acts as
an almost isometry on all $O(r)$-sparse vectors. Precisely, we
define the restricted isometry constant $\d_r$ to be the smallest
positive number such that the inequality
\begin{equation}                    \label{ri}
  C(1-\d_r) \|x\|_2^2 \le \|\Phi_T x\|_2^2 \le C(1+\d_r) \|x\|_2^2
\end{equation}
holds for some number $C>0$ and for all $x$ and all
subsets $T \subset \{1,\ldots,n\}$ of size $|T| \le r$,
where $\Phi_T$ denotes the $k \times |T|$ matrix that consists of
the columns of $\Phi$ indexed by $T$. The following theorem
is due to Candes and Tao (\cite{CT 05}, see \cite{CRTV}).

\medskip

\begin{theorem}[Restricted Isometry Condition] \label{CT}
  {\em
  Let $\Phi$ be a measurement matrix
  whose restricted isometry constant satisfies
  \begin{equation}                  \label{ric}
    \d_{3r} + 3\d_{4r} \le 2.
  \end{equation}
  Let $f$ be an $r$-sparse signal.
  Then the solution to the linear program \eqref{convex}
  is unique and is equal to $f$.
  }
\end{theorem}

\medskip

This theorem says that under the restricted isometry
condition \eqref{ric} on the measurement matrix $\Phi$,
the reconstruction problem \eqref{non-convex} is equivalent to its
convex relaxation \eqref{convex} for all $r$-sparse functions $f$.

A problem with the use of Theorem~\ref{CT} is that the restricted
isometry condition \eqref{ric} is usually difficult to check.
Indeed, the number of sets $T$ involved in this condition is exponential in $r$.
As a result, no explicit construction of a measurement matrix is presently
known that obeys the restricted isometry condition \eqref{ric}.
All known constructions of measurement matrices are randomized.

\section{Reconstruction from Fourier measurements} \label{s:Fourier}

Our goal will be to reconstruct an $r$-sparse signal $f \in \C^n$
from its discrete Fourier transform evaluated at $k = k(r,n)$
points. These points will be chosen at random and uniformly in
$\{0,\ldots,n-1\}$, forming a set $\Omega$.

The Discrete Fourier transform $\hat{f} = \Psi f$ is defined by
the DFT matrix $\Psi$ with entries
$$
\Psi_{\w,t} = \frac{1}{\sqrt{n}} \exp(-i 2 \pi \w t / n),
\ \ \ \w, t \in \{0,\ldots,n-1\}.
$$
So, our measurement matrix $\Phi$ is the submatrix of $\Psi$
consisting of random rows (with indices in $\Omega$). To be able
to apply Theorem \ref{CT}, it is enough to check that the
restricted isometry condition \eqref{ric} holds for the random
matrix $\Phi$ with high probability. The problem is -- what is the
smallest number of rows $k(r,n)$ of $\Phi$ for which this holds?
With that number, Theorem \ref{CT} immediately implies the
following reconstruction theorem for Fourier measurements:

\medskip

\begin{theorem}[Reconstruction from Fourier measurements]                             \label{Fourier rec}
  {\em
  A random set $\Omega \in \{0,\ldots,n-1\}$ of size $k(r,n)$
  satisfies the following with high probability.
  Let $f$ be an $r$-sparse signal in $\C^n$.
  Then $f$ can be exactly reconstructed from the values
  of its Fourier transform on $\Omega$ as a solution to the
  linear program
  $$
  \text{minimize }  \|f^*\|_1
  \text{ subject to } \hat{f^*}(\w) = \hat{f}(\w),
  \ \ \ \w \in \Omega.
  $$
  }
\end{theorem}

\medskip

The central remaining problem,
what is the smallest value of $k(r,n)$, is still open.
The best known estimate is due to Candes and Tao \cite{CT}:
\begin{equation}                            \label{CT k}
  k(r,n) = O(r \log^6 n).
\end{equation}
The conjectured optimal estimate would be $O(r \log n)$, 
which is known to hold for nonuniveral measuremets, i.e. for 
{\em one} sparse signal $f$ and for a random set $\Omega$
\cite{CRT}.

In this paper, we improve on the best known bound \eqref{CT k}:

\medskip

\begin{theorem}[Sample size]                \label{sample size}
  {\em
  Theorem \ref{Fourier rec} holds with
  $$
  k(r,n) = O(r \log(n) \cdot \log^2(r) \log(r \log n) ).
  $$
  }
\end{theorem}

\medskip

The dependence on $n$ is thus optimal within the $\log \log n$ factor
and the dependence on $r$ is optimal within the $\log^3 r$ factor.
So, our estimate is especially good for small $r$,
but our estimate always yields $k(r,n) = O(r \log^4 n)$.

\medskip

\remark Our results hold for transforms more general than the
discrete Fourier transform. One can replace the DFT matrix $\Psi$
by any orthogonal matrix with entries of magnitude
$O(1/\sqrt{n})$. Theorems~\ref{Fourier rec} and \ref{sample size}
hold for any such matrix.

\medskip

In the remainder of this section, we prove Theorem \ref{sample size}.
Let $\Omega$ be a random subset of $\{0,\ldots,n\}$ of size $k$.
Recall that the measurement matrix $\Phi$
that consists of the rows of $\Psi$ whose indices are in $\Omega$).
In view of Theorem \ref{ri}, it suffices to prove that the
restricted isometry constant $\d_r$ of $\Phi$ satisfies
\begin{equation}                            \label{Edelta}
  \E \d_r \le \e
\end{equation}
whenever
\begin{equation}                        \label{k}
k \ge C \Big( \frac{r \log n}{\e^2} \Big)
  \log \Big( \frac{r \log n}{\e^2} \Big)
  \log^2 r,
\end{equation}
where $\e > 0$ is arbitrary, and $C$ is some absolute
constant.

Let $y_1, \ldots, y_k$ denote the rows of the matrix $\Psi$.
Dualizing \eqref{ri} we see that \eqref{Edelta} is equivalent
to the following inequality:
$$
\E \sup_{|T| \le r}
\Big\| \id_{\C^T} - C' \sum_{i \in \Omega} y_i^T \otimes y_i^T \Big\|
\le \e
$$
with $C' = 1/\sqrt{C}$. Here and thereafter, for vectors $x,y \in
\C^n$ the tensor $x \otimes y$ is the rank-one linear operator
given by $(x \otimes y)(z) = \< x, y \> z$, where $\< \cdot \> $
is the canonical inner product on $\C^n$. The notation $x^T$
stands for the restriction of a vector $x$ on its coordinates in
the set $T$. The operator $\id_{\C^T}$ in \eqref{Esup} is the
identity on $\C^T$, and the norm is the operator norm for
operators on $\ell_2^T$.

The orthogonality of $\Psi$ can be expressed as
$\id_{\C^n} = \sum_{i=0}^{n-1} y_i \otimes y_i$.
We shall re-normalize the vectors $y_i$, letting
$
x_i = \sqrt{n} \ y_{i-1}.
$
Now we have $\|x_i\|_\infty = O(1)$ for all $i$.
The proof has now reduced to the following probabilistic
statement, which we interpret as a law of large numbers
for random operators.

\medskip

\begin{theorem}[Uniform Operator Law of Large Numbers]  \label{LLN}
  {\em 
  Let $x_1, \ldots, x_n$ be vectors in $\C^n$ with
  uniformly bounded entries: $\|x_i\|_\infty \le K$ for all $i$.
  Assume that
  $
  \id_{\C^n} = \frac{1}{n} \sum_{i=1}^n x_i \otimes x_i.
  $
  Let $\Omega$ be a random subset of $\{1,\ldots,n\}$
  of size $k$. Then
  \begin{equation}                    \label{Esup}
  \E \sup_{|T| \le r}
  \Big\| \id_{\C^T} - \frac{1}{k} \sum_{i \in \Omega} x_i^T \otimes x_i^T \Big\|
  \le \e
  \end{equation}
  provided $k$ satisfies \eqref{k} (with constant $C$ that may depend on $K$).
  }
\end{theorem}

\medskip

Theorem \ref{LLN} is proved by the techniques developed in Probability
in Banach spaces. The general roadmap is similar ton \cite{Ru}, \cite{RuFA}.
We first observe that
$$
\E \, \frac{1}{k} \sum_{i \in \Omega} x_i^T \otimes x_i^T =
\frac{1}{n} \sum_{i=1}^n x_i^T \otimes x_i^T = \id_{\C^n},
$$
so the random operator whose norm we estimate in \eqref{Esup} has
mean zero. Then the standard symmetrization (see \cite{LT} Lemma
6.3) implies that the left-hand side of \eqref{Esup} does not
exceed
$$
 2 \, \E \sup_{|T| \le r}
     \Big\| \frac{1}{k} \sum_{i \in \Omega}
             \e_i \; x_i^T \otimes x_i^T \Big\|
$$
where $(\e_i)$ are independent symmetric $\{-1,1\}$-valued random variables;
also (jointly) independent of $\Omega$.
Then the conclusion of Theorem \ref{LLN} will be easily deduced
from the following lemma.

\medskip

\begin{lemma}                   \label{lemma norm small}
  Let $x_1, \ldots, x_k$, $k \le n$, be vectors in $\C^n$ with
  uniformly bounded entries, $\|x_i\|_\infty \le K$ for all $i$.
  Then
  \begin{equation}              \label{E required}
  \E \sup_{|T| \le r}
     \Big\| \sum_{i=1}^k \e_i \; x_i^T \otimes x_i^T \Big\|
  \le k_1
      \sup_{|T| \le r} \Big\| \sum_{i=1}^k x_i^T \otimes x_i^T \Big\|^\frac{1}{2}
  \end{equation}
  where
  $
  k_1 \le C_1(K) \sqrt{r} \log(r) \sqrt{\log n} \sqrt{\log k}.
  $
\end{lemma}

\medskip

Let us show how Lemma \ref{lemma norm small} implies Theorem \ref{LLN}.
We first condition on a choice of $\Omega$ and apply Lemma \ref{lemma norm small}
for $x_i$, $i \in \Omega$. Then we take the expectation with respect to $\Omega$.
We then use the a consequence of H\"older inequality,
$\E(|X|^\frac{1}{2}) \le (\E|X|)^\frac{1}{2}$ and the triangle inequality.
Let us denote the left hand side of \eqref{Esup} by $E$. We obtain:
$$
E \le \frac{2 k_1}{\sqrt{k}}
      \E \sup_{|T| \le n}
      \norm{ \frac{1}{k} \sum_{i\in \Omega} x_i^T \otimes x_i^T}^\frac{1}{2}
  \le \frac{2 k_1}{\sqrt{k}} (E+1)^\frac{1}{2}.
$$
It follows that
$E \le C_2 \frac{2 k_1}{\sqrt{k}}$,
provided that $\frac{2 k_1}{\sqrt{k}} = O(1)$.
Theorem~\ref{LLN} now follows from our choice of $k = k(r,n)$.

\medskip

Hence it is only left to prove Lemma \ref{lemma norm small}.
Throughout the proof, $B_p^n$ and $B_p^T$
denote the unit ball of the norm $\|\cdot\|_p$ on $\C^n$.
To this end, we first replace Bernoulli r.v.'s $\e_i$ by
standard independent normal random variables $g_i$, using a
comparison principle (inequality (4.8) in \cite{LT}).
Then our problem becomes to bound the Gaussian process,
indexed by the union of the unit Euclidean balls $B_2^T$
in $\C^T$ for all subsets $I$ of $\{1,\ldots,n\}$ of size
at most $r$. We  apply Dudley's inequality (Theorem 11.17 in \cite{LT}),
which is a general upper bound on Gaussian processes.
Let us denote the left hand side of \eqref{Esup} by $E_1$. We obtain:
\begin{align*}
E_1 &\le 
   C_3 \E \sup_{|T| \le r}
     \Big\| \sum_{i=1}^k g_i \; x_i^T \otimes x_i^T \Big\| \\
&= C_3 \E \sup_{\substack{|T| \le r \\ x \in B_2^T}}
    \Big| \sum_{i=1}^k g_i \< x_i, x\> ^2 \Big| \\
&\le C_4 \int_0^\infty \log^{1/2} N \big( \cup_{|T| \le r} B_2^T, \d, u \big) \; du,
\end{align*}
where $N(Z,\d,u)$ denotes the minimal number of balls of radius $u$
in metric $\d$ centered in points of $Z$, needed to cover the set $Z$.
The metric $\d$ in Dudley's inequality
is defined by the Gaussian process, and in our case it is
\begin{align*}
\d(x,y)
&= \Big[ \sum_{i=1}^M \big( \< x_i, x\> ^2 - \< x_i, y\> ^2 \big)^2
   \Big]^\frac{1}{2} \\
&\le \Big[ \sum_{i=1}^k \big( \< x_i, x\> + \< x_i, y\> \big)^2
   \Big]^\frac{1}{2}
   \ \max_{i \le k} |\< x_i, x-y \> | \\
&\le 2 \max_{\substack{|T| \le r \\ z \in B_2^T}}
  \Big[ \sum_{i=1}^k \< x_i, z\> ^2 \Big]^\frac{1}{2}
  \ \max_{i \le k} |\< x_i, x-y \> | \\
&= 2 R \max_{i \le k} |\< x_i, x-y \> |,
\end{align*}
where
$$
R := \sup_{|T| \le r} \Big\| \sum_{i=1}^k x_i^T \otimes x_i^T \Big\|^\frac{1}{2}.
$$
Hence
\begin{equation}                \label{with D2}
E_1 \le C_5 R \sqrt{r}
  \int_0^\infty \log^{1/2}
    N \big( \frac{1}{\sqrt{r}} D_2^{r,n}, \|\cdot\|_X, u \big) \; du.
\end{equation}
Here
$$
D_p^{r,n} = \bigcup_{|T| \le r}  B_p^T, \ \ \ \|x\|_X = \max_{i
\le k} |\< x_i, x\> |.
$$
We will use containments
\begin{equation}                \label{D containments}
\frac{1}{\sqrt{r}} D_2^{r,n} \subseteq D_1^{r,n} \subseteq K B_X,
\ \ \ D_1^{r,n} \subseteq B_1^n,
\end{equation}
where $B_X$ denotes the unit ball of the norm $\|\cdot\|_X$. The
second containment follows from the uniform boundedness of
$(x_i)$. We can thus replace $\frac{1}{\sqrt{r}} D_2^{r,n}$ in
\eqref{with D2} by $D_1^{r,n}$. Comparing \eqref{with D2} to the right hand side
of \eqref{E required} we see that, in order to complete the proof
of Lemma \ref{lemma norm small}, it suffices to show that
\begin{equation}                            \label{integral required}
\int_0^K \log^{1/2}
  N \big( D_1^{r,n},  \|\cdot\|_X, u \big) du
\le C_6 \log(r) \sqrt{\log n} \sqrt{\log k},
\end{equation}
with $C_6 = C_6(K)$.
To this end, we will estimate the covering numbers in this integral
in two different ways. For big $u$, we will just use the second
containment in \eqref{D containments}, which allows us to
replace $D_1^{r,n}$ by $B_1^n$.

\medskip

\begin{lemma}   \label{lemma covering numbers}
  Let $x_1, \ldots, x_k$, $k \le n$, be vectors as in Lemma \ref{lemma norm small}.
  Then for all $u > 0$ we have
  \[
  N(B_1^n, \|\cdot\|_X, u) \le (2n)^m,
  \]
  where $m = C_7 K^2 \log(k)/u^2$.
\end{lemma}

\medskip

\begin{proof}
We use the empirical method of Maurey. Fix a vector $y \in B_1^n$.
Define a random vector $Z \in \R^n$ that takes values
$(0,\ldots,0,\sign(y(i)),0,\ldots,0)$ with probability $|y(i)|$
each, $i=1,\ldots,n$ (all entries of that vector are zero except
$i$-th). Here $\sign(z)=z/|z|$, whenever $z \neq 0$, and $0$
otherwise. Note that $\E Z = y$. Let $Z_1,\ldots,Z_m$ be
independent copies of $Z$. Using symmetrization as before, we see
that
\[
E_3 := \E \Big\| y - \frac{1}{m} \sum_{j=1}^m Z_j \Big\|_X
\le \frac{2}{m} \ \E \Big\| \sum_{j=1}^m \e_j Z_j \Big\|_X.
\]
Now we condition on a choice of $(Z_j)$ and take the
expectation with respect to random signs $(\e_j)$.
Using comparison to Gaussian variables as before, we obtain
\begin{align*}
E_4
  &:=  \E \Big\| \sum_{j=1}^m \e_j Z_j \Big\|_X
  \le C_7 \E \Big\| \sum_{j=1}^m g_j Z_j \Big\|_X \\
  &=    C_7 \E \max_{i \le k} \Big| \sum_{j=1}^m g_j \pr{Z_j}{x_i} \Big|.
\end{align*}
For each $i$, $\gamma_i := \sum_{j=1}^m g_j \pr{Z_j}{x_i}$
is a Gaussian random variable with zero mean and with variance
$$
\sigma_i = \big( \sum_{j=1}^m |\pr{Z_j}{x_i}|^2 \big)^{1/2} \le K
\sqrt{m},
$$
since $|\pr{Z_j}{x_i}| \le \|x_i\|_{\infty} \le K$. Using a simple
bound on the maximum of Gaussian random variables (see (3.13) in
\cite{LT}), we obtain
$$
E_4 \le C_7 \E \max_{i \le k} |\gamma_i| \le C_8 \sqrt{\log k}
\max_{i \le k} \sigma_i \le C_8 \sqrt{\log k} K \sqrt{m}.
$$
Taking the expectation with respect to $(Z_j)$ we obtain
$$
E_3 \le \frac{2}{m} \E(E_4) \le \frac{2C_8 K \sqrt{\log
k}}{\sqrt{m}}.
$$
With the choice of $m$ made in the statement of the lemma,
we conclude that
$
E_3 \le u.
$
We have shown that for every $y \in B_1^n$, there exists a $z \in
\C^n$ of the form $z = \frac{1}{m} \sum_{j=1}^m Z_j$ such that
$\|y-z\|_X \le u$. Each $Z_j$ takes $2n$ values, so $z$ takes
$(2n)^m$ values. Hence $B_1^n$ can be covered by at $(2n)^m$ balls
of norm $\|\cdot\|_X$ of radius $u$. A standard argument shows
that we can assume that these balls are centered in points of
$B_1^n$. This completes the proof of Lemma~\ref{lemma covering
numbers}.
\end{proof}

\medskip

For small $u$, we will use a simple volumetric estimate.
The diameter of $B_1^r$ considered as a set in $\C^n$ is
at most $K$ with respect to the norm $\|\cdot\|_X$
(this was stated as the last containment in \eqref{D containments}).
It follows that $N(B_1^r,\|\cdot\|,u) \le (1+2K/u)^r$
for all $r > 0$, see (5.7) in \cite{Pi}.
The set $D_1^{r,n}$ consists of $d(r,n) = \sum_{j=1}^r \binom{n}{i}$
balls of form $B_1^T$, thus
\begin{equation}                    \label{volumetric}
N \big( D_1^{r,n},  \|\cdot\|_X, u \big)
\le d(n,r) (1+2K/u)^r.
\end{equation}

Now we combine the estimate of the covering number $N(u) =
\log^{1/2} N \big( D_1^{r,n},  \|\cdot\|_X, u \big)$ of Lemma 3.6,
and the volumetric estimate \eqref{volumetric}, to bound the
integral in \eqref{integral required}. Using Stirling's
approximation, we see that $d(r,n) \le (C_{9} n/r)^r$. Thus
\begin{gather*}
N(u) \le C_{10} \sqrt{r} \big[ \sqrt{\log(n/r)} + \sqrt{\log(1+2/u)} \big] =: N_1(u), \\
N(u) \le \frac{C_{10}}{u} \sqrt{\log k} \sqrt{\log n} =: N_2(u),
\end{gather*}
where $C_{10} = C_{10}(K)$.
Then we bound the integral in \eqref{integral required} as
\begin{align*}
\int_0^K N(u) \; du
&\le \int_0^A N_1(u) \; du + \int_A^K N_2(u) \; du \\
&\le C_{11} A \sqrt{r} \big[ \sqrt{\log(n/r)} + \log(1+2/A) \big] \\
  &\ \ + C_{11} \log(1/A) \sqrt{\log k} \sqrt{\log n},
\end{align*}
where $C_{11} = C_{11}(K)$. Choosing $A = 1/\sqrt{r}$, we conclude
that the integral in \eqref{integral required} is at most
$
\sqrt{\log(n/r)} + \log r + \log(r) \sqrt{\log k} \sqrt{\log n}.
$
This proves \eqref{integral required}, which completes the proof
of Lemma \ref{lemma norm small} and thus of Theorems~\ref{LLN}
and \ref{sample size}.
\endproof

\section{Reconstruction from Gaussian measurements}  \label{s:Gauss}

Our goal will be to reconstruct an $r$-sparse signal $f \in \R^n$
from $k = k(r,n)$ Gaussian measurements. These are given by
$\Phi f \in \R^k$, where $\Phi$ is a $k \times n$ random matrix
(`Gaussian matrix' in the sequel),
whose entries are independent $N(0,1)$ random variables.
The reconstruction will be achieved by solving the linear
program \eqref{convex}.

The problem again is to find the smallest number of measurements
$k(r,n)$ for which, with high probability, we have an exact
reconstruciton of every $r$-sparse signal $f$ from its
measurements $\Phi f$? It has recently been shown in
\cite{CT 05, RV, CRTV} that
\begin{equation}                    \label{old k}
k(r,n) = O(r \log(n/r)),
\end{equation}
and was extended in \cite{MPT} to sub-gaussian measurements.
This is asymptotically optimal. However, the constant factor
implicit in \eqref{old k} has not been known; previous proofs of
\eqref{old k} yield unreasonably weak constants (of order $2,000$
and higher). In fact,
{\em there has not been known any
  theoretical guarantees with reasonable constants for Linear
  Programming based reconstructions}.
So, there is presently a gap
between theoretical guarantees and good practical performance of
reconstruction \eqref{convex} (see e.g. \cite{CRTV}). Here we shall prove
a first practically reasonable guarantee of the form \eqref{old
k}:
\begin{align}                       \label{new k}
k(r,n) \le c_1 r \big[ c_2 + \log(n/r) \big] (1 + o(1)),  \\
c_1 = 6 + 4 \sqrt{2} \approx 11.66,
\ \ c_2 = 1.5. \notag
\end{align}

\medskip

\begin{theorem}[Reconstruction from Gaussian measurements]                      \label{Gaussian rec}
  {\em
  A $k \times n$ Gaussian matrix $\Phi$ with $k > k(r,n)$
  satisfies the following with probability
  $$
  1 - 3.5 \; \exp \Big( -\big( \sqrt{k}-\sqrt{k(r,n)} \big)^2/18 \Big).
  $$
  Let $f$ be an $r$-sparse signal in $\R^n$.
  Then $f$ can be exactly reconstructed from the measurements $\Phi f$
  as a unique solution to the linear program \eqref{convex}.
  }
\end{theorem}

\medskip

Our proof of Theorem \ref{Gaussian rec} is direct, we will not use
the Restricted Isometry Theorem \ref{CT}. 
The first part of this argument follows a general method of \cite{MPT}. 
One interprets the exact reconstruction
as the fact that the (random) kernel of $\Phi$ misses the cone generated by 
the (shifted) ball of $\ell_1$. Then one embeds the cone in a universal 
set $D$, which is easier to handle, and proves that the random subspace 
does not intersect $D$. However, to obtain good constants as in \eqref{new k},
we will need to (a) improve the constant of embedding into $D$ from \cite{MPT}, 
and (b) use Gordon's Escape Through the Mesh Theorem \cite{G}, 
which is tight in terms of constants.
In Gordon's theorem, one measures
the size of a set $S$ in $\R^n$ by its {\em Gaussian width}
$$
w(D) = \E \sup_{x \in S} \pr{g}{x},
$$
where $g$ is a random vector in $\R^n$ whose components
are independent $N(0,1)$ random variables (Gaussian vector).
The following is Gordon's theorem \cite{G}.

\medskip

\begin{theorem}[Escape Through the Mesh (Gordon)]   \label{escape}
  {\em
  Let $S$ be a subset of the unit Euclidean sphere $S^{n-1}$ in $\R^n$.
  Let $Y$ be a random $(n-k)$-dimensional subspace of $\R^n$,
  distributed uniformly in the Grassmanian with respect to the Haar measure.
  Assume that
  $
    w(S) > \sqrt{k}.
  $
  Then
  $
  Y \cap S = \emptyset
  $
  with probability at least
  $$
  1 - 3.5 \; \exp \Big( -\big( k/\sqrt{k+1} - w(S) \big)^2/18 \Big).
  $$
  }
\end{theorem}

\medskip

We will now prove Theorem \ref{Gaussian rec}.
First note that the function $f$ is the unique solution
of \eqref{convex} if and only if $0$ is the unique solution of the problem
\begin{equation}                    \label{convex'}
  \text{minimize }  \|f-g^*\|_1 \text{ subject to } \Phi g^* \in \Ker(\Phi) =: Y.
\end{equation}
$Y$ is a $(n-k)$-dimensional subspace of $\R^n$. Due to the rotation invariance
of the Gaussian random vectors, $Y$ is distributed uniformly in the Grassmanian
$G_{n-k,n}$ of $(n-k)$-dimensional subspaces of $\R^n$, with respect to the
Haar measure.

Now, $0$ is the unique solution to \eqref{convex'}
if and only if $0$ is the unique metric projection of $f$
onto the subspace $Y$ in the norm $\|\cdot\|_1$.
This in turn is equivalent to the fact that $0$ is the unique
contact point between the subspace $Y$ and the ball of the norm $\|\cdot\|_1$
centered at $f$:
\begin{equation}                    \label{convex''}
(f + \|f\|_1 B_1^n) \cap Y = \{0\}.
\end{equation}
(Recall that $B_p^n$ is the unit ball of the norm $\|\cdot\|_p$.)
Let $\mathcal{C}_f$ be the cone in $\R^n$ generated by the set $f + \|f\|_1 B_1^n$
(the cone of a set $A \in \R^n$ is defined as $\{ta \mid a \in A, \ t \in \R^+\}$).
Then the statement that \eqref{convex''} holds for all $r$-sparse functions $f$
is clearly equivalent to
\begin{equation}                        \label{miss cone}
  \mathcal{C}_f \cap Y = \{0\}
  \ \ \text{for all $r$-sparse functions $f$}.
\end{equation}
We can represent the cone $\mathcal{C}_f$ as follows.
Let
$$
T^+ = \{i \mid f(i)>0\}, \
T^- = \{j \mid f(i)<0\}, \
T   = T^+ \cup T^-.
$$
Then
\[
  \mathcal{C}_f = \Big\{t \in \R^n
    \mid \sum_{i \in T^-} t(i) - \sum_{i \in T^+} t(i)
      + \sum_{i \in T^c} |t(i)| \le 0
             \Big\}.
\]
We will now bound the cone $\mathcal{C}_f$ by a universal set, which
does not depend on $f$.

\medskip

\begin{lemma} \label{l: inclusion}
  Consider the spherical part of the cone,
  $K_f = \mathcal{C}_f \cap S^{n-1}$.
  Then $K_f \subset (\sqrt{2}+1) D$, where
  $$
  D = \conv \{x \in S^{n-1} \mid |\supp(x)| \le r\}.
  $$
\end{lemma}

\medskip

\begin{proof}
Fix a point $x \in \mathcal{C} \cap S^{n-1}$.
We have
$$
\sum_{i \in T}|x(i)| \le \sqrt{|I|} \le \sqrt{r}, \ \
\sum_{i \in T^c}|x(i)| \le \sum_{i \in T}|x(i)| \le \sqrt{r}.
$$
The norm $\|\cdot\|_D$ on $\R^n$ whose unit ball is $D$
can be computed as
\[
  \|x\|_D=
  \sum_{l=1}^L \Big(\sum_{i \in I_l} (x(i)^*)^2 \Big)^{1/2},
\]
where $L=\lceil n/r \rceil$, $I_l=\{r(l-1)+1 \etc rl\}$, for
$l<L$, $I_L=\{r(L-1)+1 \etc n\}$, and $(x(i)^*)$ is a
non-decreasing rearrangement of the sequence $(|x(i)|)$.

Set $F = F(x) = \{i \mid |x(i)| \ge 1/\sqrt{r} \}$.
Since $x \in S^{n-1}$, we have $|F| \le r$.
Hence, for any $x \in K$ there exists a set
$E = E(x) \subset \{1 \etc m \}$, which consists of $2r$ elements
and such that $E \supseteq F \cup I$.
Therefore, $x$ can be represented as $x=x'+x''$
so that  $\supp(x') \subseteq E$,
$\|x\|_2 \le 1$,
$\supp (x'') \subseteq E^c$,
$\|x''\|_\infty \le 1/\sqrt{r}$.
Set
\[
   V_{E}=B_2^{E} \times \Big (\sqrt{r}B_1^{E^c} \cap
  \frac{1}{\sqrt{r}} B_{\infty}^{E^c} \Big ).
\]
Then the above argument shows that
$
   K_f \subset \bigcup_{|E|=2r} V_E=:W.
$

The maximum of $\|x\|_D$ over $x \in W$ is attained at the
extreme points of the sets $V_E$, which have the form
$x=x'+x''$, where $x' \in S^{E}$, and $x''$ has coordinates 0 and
$\pm 1/\sqrt{r}$ with $r$ non-zero coordinates. Notice that since
$|\supp (x')| \le 2r$, $\|x'\|_D \le \sqrt{2} \|x'\|_2$.
Thus, for any extreme point $x$ of $V_{E}$,
\[
  \|x\|_D \le  \|x'\|_D+ \|x''\|_D
  \le \sqrt{2} \|x'\|_2+\|x''\|_2
  \le \sqrt{2}+1.
\]
The second inequality follows from $\supp(x') \le 2r$ and
$\supp(x'')=r$. This completes the proof of the lemma.
\end{proof}

To use Gordon's escape through the mesh theorem,
we have to estimate the Gaussian width of $D$.

\medskip

\begin{lemma} \label{l: mean width}
\[
   w(D) \le \sqrt{2r \log(e^{3/2}n / r)} (1+o(1)).
\]
\end{lemma}

\medskip

\begin{proof}
By definition,
\[
w(D) = \sup_{|J|=r} \Big( \sum_{i \in J} |g(i)|^2 \Big)^{1/2}.
\]
Let $p > 1$ be a number to be chosen later.
By H\"older's inequality, we have
\begin{align*}
  w(D)
  &\le \E \Big( \sum_{|J|=r} \Big( \sum_{i \in J} |g(i)|^2
            \Big)^{p/2} \Big)^{1/p} \\
  &\le \binom{n}{r}^{1/p} \Big( \E \Big( \sum_{i=1}^r |g(i)|^2
            \Big)^{p/2} \Big)^{1/p} \\
  &\le \Big( \frac{en}{r} \Big)^{r/p}
     \Big( 2^{p/2} \cdot \frac{\Gamma (p/2+r/2)}{\Gamma(r/2)}
     \Big)^{1/p}.
\end{align*}
By the Stirling's formula,
\[
    2^{p/2} \cdot \frac{\Gamma (p/2+r/2)}{\Gamma(r/2)}
    = \left (1+ \frac{p}{r} \right )^{\frac{r+1}{2}}
    \left (\frac{p+r}{e} \right )^{p/2}
    (1+o(1)).
\]
Therefore,
$
  w(D) \le \left ( \frac{en}{r} \right )^{r/p}
    \left ( \frac{p+r}{e} \right )^{1/2}
    (1+o(1)).
$
Now set $p=2r \log (\frac{en}{r} )$. Then
\[
  w(D) \le (p+r)^{1/2} (1+o(1))
  = \sqrt{2r \log \frac{e^{3/2}n}{r}} (1+o(1)).
\]
\end{proof}

To deduce \eqref{miss cone}
we define
$
S = \bigcup_f K_f,
$
where the union is over all $r$-sparse functions $f$.
Then \eqref{miss cone} is equivalent to
\begin{equation}                        \label{miss S}
S \cap Y = \emptyset.
\end{equation}
Lemma~\ref{l: inclusion} implies that
$S \subseteq (\sqrt{2} + 1)D$.
Then by Lemma~\ref{l: mean width},
\begin{align*}
w(S)
\le (\sqrt{2}+1) w(D) 
= (1-o(1)) \sqrt{k(r,n)}.
\end{align*}
Then \eqref{miss S} follows Gordon's Theorem~\ref{escape}.
This completes the proof of Theorem~\ref{Gaussian rec}.
\endproof

{\bf Acknowledgement.} After this paper was announced, A.Pajor 
pointed out that Lemma 3.6 was proved by B.Carl in \cite{Carl},
see Prop.3 and below. We also thank Emmanuel Candes for important
remarks.

\end{document}